\documentclass{article}

\usepackage[square,sort,numbers]{natbib}
\usepackage[utf8]{inputenc} 
\usepackage[T1]{fontenc}    
\usepackage{hyperref}       
\usepackage{url}            
\usepackage{booktabs}       
\usepackage{amsfonts}       
\usepackage{nicefrac}       
\usepackage{microtype}      
\usepackage{amssymb}
\usepackage{amsmath}
\usepackage{amsthm}
\usepackage{graphicx}

\title{Tutorial: Complexity analysis of Singular Value Decomposition and its variants}

\author{%
  Xiaocan Li \space\space\space Shuo Wang \space\space\space Yinghao Cai\\
  Institute of Automation, Chinese Academy of Sciences\\
  University of Chinese Academy of Sciences\\
  \texttt{\{lixiaocan2017, shuo.wang, yinghao.cai\}@ia.ac.cn} \\
}

\begin{document}

\maketitle

\begin{abstract}

We compared the regular Singular Value Decomposition (SVD), truncated SVD, Krylov method and Randomized PCA, in terms of time and space complexity. It is well-known that Krylov method and Randomized PCA only performs well when $k \ll n$, i.e. the number of eigenpair needed is far less than that of matrix size. We compared them for calculating all the eigenpairs. We also discussed the relationship between Principal Component Analysis and SVD.
\end{abstract}

\section{Introduction}
Dimensionality reduction has always been a trendy topic in machine learning. Linear subspace method for reduction, e.g., Principal Component Analysis and its variation have been widely studied\cite{Xu1995Robust, Partridge2002Robust, Zou2006Sparse}, and some pieces of literature introduce probability and randomness to realize PCA\cite{Tipping2010Probabilistic, Halko2010Finding, Halko2011An, J1989Estimating, Woolfe2008A}. However, linear subspace is not applicable when the data lies in a non-linear manifold\cite{Belkin2014Laplacian, Kai2010Clustered}. Due to the direct connection with PCA, Singular Value Decomposition (SVD) is one of the most well-known algorithms for low-rank approximation\cite{Woodruff2014Low, Frieze1998Fast}, and it has been widely used throughout the machine learning and statistics community. 
Some implementations of SVD are solving least squares\cite{Golub1970Singular, Zhang2010Regularized}, latent semantic analysis\cite{deerwester1990indexing, hofmann2001unsupervised}, genetic analysis, matrix completion\cite{Recht2012Exact, Cai2008A, Candes2010Matrix, Cand2010The}, data mining\cite{Chong2016Feature, Belabbas2009On} etc. However, when it comes to a large scale matrix, the runtime of traditional SVD is intolerable and the memory usage could be enormously consuming. 

$\textbf{Notations}$: we have a matrix $A$ with size $m\times n$, usually $m \gg n$. Our goal is to find the Principal Components (PCs) given the cumulative explained variance threshold $t$.

$\textbf{Assumptions}$: In this tutorial, every entry of matrix $A$ is real-valued; W.l.o.g., assume $m \gg n$ and $A$ has zero mean over each feature. 

For your information, either each column or row of $A$ could represent an example, and the definition will be specified when necessary.

In the traditional approach of PCA, we need to compute the covariance matrix $S$, then perform the eigen-decomposition on $S$. By selecting the top $K$ largest eigenvalues and corresponding eigenvectors, we get our Principal Components (PCs). Nevertheless, if each column of $A$ is an example and the row size of $A$ is tremendously large, saving even larger covariance matrix into memory is expensive, let alone the eigen-decomposition process.

\section{Preliminary knowledge}\label{sec:pre}
\subsection{Singular Value Decomposition (SVD)}
Any real or complex matrix can be approximated over the summation of a series of rank-1 matrix. In SVD, we have
\begin{align*}
A = U\Sigma V^T
\end{align*}
where
\begin{align*}
U &= [u_1 \cdots u_n]\in\mathbb{R}^{m \times n}\\
\Sigma &= diag(\sigma_1 \cdots \sigma_n)\in\mathbb{R}^{n\times n}\\
V &= [v_1 \cdots v_n]\in\mathbb{R}^{n \times n}
\end{align*}
Here $U$ and $V$ are orthogonal matrices, i.e.
\begin{align*}
U^T U &= I_n\\
UU^T &= I_m\\
V^TV&=VV^T=I_n
\end{align*}
We could also rewrite SVD as following
\begin{align}\label{AV}
Av_i = \sigma_i u_i, i = 1 \cdots n 
\end{align}
Geometrically speaking, the matrix $A$ rotates the unit vector $v_i$ to $u_i$ and then stretches the Euclidean norm of $u_i$ with a factor of $\sigma_i$.

The orthogonal matrices $U$ and $V$ can be obtained by eigen-decomposition of matrix $AA^T$ and $A^TA$, and the singular values $\sigma_i's$ are the square root of the eigenvalues of $AA^T$ or $AA^T$.

Proof:
\begin{align} 
AA^T &= U\Sigma V^T V\Sigma^T U^T = U\Sigma^2 U^T \label{AA}\\
A^TA &= V\Sigma^T U^T U\Sigma V^T = V\Sigma^2 V^T \label{AB}
\end{align}

Set $\Lambda = \Sigma^2= diag(\sigma_1^2 \cdots \sigma_n^2)$, i.e. 

\begin{align}\label{LambdaSigma}
\lambda_i = \sigma_i^2
\end{align}

We could rewrite Eq. (\ref{AA}) and Eq. (\ref{AB}) as

\begin{align} 
AA^Tu_i &= \lambda_i u_i \\
A^TAv_i &= \lambda_i v_i
\end{align}

Therefore, the column vectors of $U$ and $V$ are the eigenvectors (with the unit norm) of $AA^T$ and $A^T A$, respectively. Moreover, the eigenvalues $\lambda_i$ are square of singular values $\sigma_i$, as in Eq. (\ref{LambdaSigma}). In other words, the square root of eigenvalues are singular values. 
\subsection{Relations with PCA}
\textbf{If each column of $A$ represents an example or data point}, set $Y=U^TA$ to our transformed data points, where $U$ is the left singular matrix in SVD. The covariance matrix of $Y$ is
\begin{align}
cov(Y) = YY^T = U^TAA^TU = U^TU\Sigma^2U^TU = \Sigma^2
\end{align}

\textbf{If each row of $A$ represents an example or data point}, set $Y=AV$ to our transformed data points, where $V$ is the right singular matrix in SVD. The covariance matrix of $Y$ is
\begin{align}
cov(Y) = Y^TY = V^TA^TAV = V^TV\Sigma^2V^TV = \Sigma^2
\end{align}
It means the transformed data $Y$ are uncorrelated. Therefore, the column vectors of orthogonal matrix $U$ in SVD are the projection bases for Principal Components.

\subsection{Truncated SVD}
Although the derivation of SVD is clear theoretically, practically speaking, however, it is unwise to do eigen-decomposition on matrix $AA^T$, as it has a tremendous size of $m \times m$, which will deplete memory and cost a great amount of time. On the contrary, the matrix $A^TA$ only has a size of $n \times n$, thus it is plausible that we compute orthogonal matrix $V$ first.
Then we can plug the Eq. (\ref{AV}) into Eq. (\ref{LambdaSigma}) and get
\begin{align}
u_i = \dfrac{Av_i}{\sqrt{\lambda_i}}, i = 1 \cdots n
\end{align}

This equation is the key to improving time and space efficiency because we do not perform eigen-decomposition on huge matrix $AA^T\in\mathbb{R}^{m\times m}$, which takes $\mathcal{O}(m^3)$ time and $\mathcal{O}(m^2)$ space.

Then we column-wisely combine $u_i$ to get $U$, and the same for $v_i$ to form $V$. For $\Sigma$, it is $diag(\sqrt{\sigma_1} \cdots \sqrt{\sigma_n})$.

\subsection{PCA Evaluation}
In Section \ref{sec:pre}, we have proved that the column vectors of orthogonal matrix $U$ or $V$ in SVD are the projection bases for PCA. In the literature of PCA, there are many criteria for evaluating the residual error, e.g., Frobenius norm and induced $L_2$ norm of the difference matrix (original matrix minus approximated matrix), explained variance and cumulative explained variance. 

Usually, we use the cumulative explained variance criterion for evaluation.

\textbf{cumulative explained variance criterion:} Given the threshold $t$, find the minimal integer $K$ such that
\begin{align}
\frac{\Sigma_{i=1}^{K} \lambda_i }{\Sigma_{i=1}^n \lambda_i}\geq t
\end{align}
where each $\lambda_i$ is the eigenvalue of matrix $A^TA$ or $AA^T$. Every $\lambda_i=\sigma_i^2$ indicates the variance in principal axis, this is why the criterion is named cumulative explained variance.

For those SVD or PCA algorithms who do not obtain all the eigenvalues or can not get accurate them accurately, it seems that the denominator term $\Sigma_{i=1}^n \lambda_i$ can not be calculated. Actually, the sum of all eigenvalues can be done by
\begin{align}\label{sumeigen}
\Sigma_{i=1}^n \lambda_i = tr(A^TA) = tr(AA^T) = \|A\|_F^2 = \Sigma_{i, j} a_{ij}^2
\end{align}
\textbf{Therefore, we do not need to implement eigen-decomposition on either large matrix $AA^T\in\mathbb{R}^{m\times m}$ or small matrix $A^TA\in\mathbb{R}^{n\times n}$.} Eq. (\ref{sumeigen}) saves us $\mathcal{O}(n^3)$ time and $\mathcal{O}(n^2)$ space.

\section{Complexity Analysis}
In this section, we compare the time complexity and space complexity of Krylov method, Randomized PCA and truncated SVD. Due to the copyrights issue, the mechanisms of MatLab and Python built-in econSVD are not available, whose complexity analysis will not be conducted.

To restate again, our matrix $A$ has size $m\times n$ and $m>n$.

\subsection{Time Complexity}
For time complexity, we use the number of FLoating-point OPerations (FLOP) as a quantification metric. 

For matrix $A \in \mathbb{R}^{m\times n}$ and $B \in \mathbb{R}^{n\times l}$, the time complexity of matrix multiplication $AB$ takes $mnl$ FLOP of products and $ml(n-1)$ FLOP of summations. Therefore, the multiplication of two matrices takes $\mathcal{O}(mnl+ml(n-1))=\mathcal{O}(2mnl-ml) = \mathcal{O}(mnl)$ FLOP. We could ignore the coefficient here for it will not bring bias to our analysis.

For your information, the coefficient of time complexity in Big-O notation will not be ignored when comparing different SVD algorithms as it is of importance in our analysis.
\subsubsection{Krylov method}
For the Krylov method, we discuss the time complexity of each step.
\begin{enumerate}
	\item Forming standard normal distribution matrix $G$ of size $n \times l$ takes $\mathcal{O}(nl)$ FLOP, where in practice $l=0.5n$.
	\item Forming matrix $H^{(0)}=AG \in \mathbb{R}^{m\times l}$ takes $\mathcal{O}(mnl)$ FLOP.
	\item Forming matrix $H^{(i)}=A(A^T H^{(i-1)})\in \mathbb{R}^{m\times l}$ takes $\mathcal{O}(2imnl)$ FLOP, as the matrix multiplication in bracket $A^T H^{(i-1)}$ takes $\mathcal{O}(nml)$ FLOP, then multiplying by $A$ takes $\mathcal{O}(mnl)$ FLOP. In total, it takes $\mathcal{O}{(2imnl)}$ FLOP to generate $i$ matrices. 
	\item Forming matrix $H=\left(H^{(0)}\left|H^{(1)}\right| \ldots\left|H^{(i-1)}\right| H^{(i)}\right) \in \mathbb{R}^{m\times (i+1)l}$ by concatenating each $H^{(i)}$ takes $\mathcal{O}(1)$ FLOP.
	\item Performing QR decomposition on $H$ takes $\mathcal{O}(2m[(i+1)l]^2 - 2[(i+1)l]^3/3)$ FLOP.
	\item Forming $T = A^T Q \in \mathbb{R}^{n\times (i+1)l}$ takes $\mathcal{O}((i+1)mnl)$ FLOP.
	\item Performing SVD on $T=\tilde{V}\tilde{\Sigma}W^T$ takes $\mathcal{O}(n[(i+1)l]^2)$ FLOP.
	\item Forming $\tilde{U}=QW$ takes $\mathcal{O}(m[(i+1)l]^2)$ FLOP.
\end{enumerate}

In total, the time complexity of Krylov method is 
\begin{align}
\mathcal{O}(nl+(3i+2)mnl+(i+1)^2 l^2 (m^2+n^2+2m-\dfrac{2}{3}(i+1)l)) \text{ FLOP}
\end{align}

In practice, $l = 0.5n, i=1$, then the time complexity will be
\begin{align}
\mathcal{O}(\dfrac{n^2}{2}+\dfrac{5}{2}mn^2+4 (m^2 n^2 +2mn^2+n^4 -\dfrac{2}{3}n^3)) \text{ FLOP}
\end{align}

\subsubsection{Randomized PCA}
For Randomized PCA, we discuss the time complexity of each step. It is very similar to the Krylov method.
\begin{enumerate}
	\item Forming standard normal distribution matrix $G$ of size $n \times l$ takes $\mathcal{O}(nl)$ FLOP, where in practice $l=0.5n$.
	\item Forming matrix $H^{(0)}=AG \in \mathbb{R}^{m\times l}$ takes $\mathcal{O}(mnl)$ FLOP.
	\item Forming matrix $H^{(i)}=A(A^T H^{(i-1)})\in \mathbb{R}^{m\times l}$ takes $\mathcal{O}(2imnl)$ FLOP, as the matrix multiplication in bracket $A^T H^{(i-1)}$ takes $\mathcal{O}(nml)$ FLOP, then multiplying by $A$ takes $\mathcal{O}(mnl)$ FLOP. In total, it takes $\mathcal{O}{(2imnl)}$ FLOP. 
	\item Performing QR decomposition on $H$ takes $\mathcal{O}(2ml^2 - 2l^3/3)$ FLOP.
	\item Forming $T = A^T Q \in \mathbb{R}^{n\times l}$ takes $\mathcal{O}(mnl)$ FLOP.
	\item Performing SVD on $T=\tilde{V}\tilde{\Sigma}W^T$ takes $\mathcal{O}(nl^2)$ FLOP.
	\item Forming $\tilde{U}=QW$ takes $\mathcal{O}(ml^2)$ FLOP.
\end{enumerate}

In total, the time complexity of Randomized PCA is 
\begin{align}
\mathcal{O}(nl+(2i+2)mnl+l^2 (3m+n-\dfrac{2}{3}l)) \text{ FLOP}
\end{align}

In practice, $l = 0.5n, i = 1$, then the time complexity will be
\begin{align}
\mathcal{O}(\dfrac{n^2}{2}+\dfrac{11}{4}mn^2+\dfrac{n^3}{6}) \text{ FLOP}
\end{align}

\subsubsection{Truncated SVD}
We discuss the time complexity of Truncated SVD for each step.
\begin{enumerate}
	\item Forming matrix $A^T A\in\mathbb{R}^{n\times n}$ takes $\mathcal{O}(mn^2)$ FLOP.
	\item Performing eigen-decomposition on $A^TA \in\mathbb{R}^{n\times n}$ takes $\mathcal{O}(n^3)$ FLOP.
	\item Taking the square root of each eigenvalue of $A^T A$ takes $\mathcal{O}(n)$ FLOP.
	\item Forming $u_i = \dfrac{Av_i}{\sigma_i}$ takes $\mathcal{O}(n(mn+m))$ FLOP, as $Av_i$ takes $\mathcal{O}(mn)$ FLOP while divided by $\sigma_i$ takes $\mathcal{O}(m)$ FLOP. In total, we have $n$ equations like this, thus it takes $\mathcal{O}(n(mn+m))$ FLOP.
\end{enumerate}

In total, the time complexity of Truncated SVD is
\begin{align}
\mathcal{O}(2mn^2+n^3+n+mn)
\end{align}

\subsection{Space Complexity}
We evaluate the space complexity by the number of matrix entries. For a matrix $A\in\mathbb{R}^{m\times n}$, its space complexity is $\mathcal{O}(mn)$. In MatLab or Python programming language, each entry takes 8 bytes memory.
\subsubsection{Krylov method}
\begin{enumerate}
	\item Forming standard normal distribution matrix $G$ of size $n \times l$ takes $\mathcal{O}(nl)$, where in practice $l=0.5n$.
	\item Forming matrix $H^{(0)}=AG \in \mathbb{R}^{m\times l}$ takes $\mathcal{O}(ml)$.
	\item Forming matrix $H^{(i)}=A(A^T H^{(i-1)})\in \mathbb{R}^{m\times l}$ takes $\mathcal{O}(ml)$. 
	\item Forming matrix $H=\left(H^{(0)}\left|H^{(1)}\right| \ldots\left|H^{(i-1)}\right| H^{(i)}\right) \in \mathbb{R}^{m\times (i+1)l}$ by concatenating each $H^{(i)}$ takes $\mathcal{O}((i+1)ml)$.
	\item Performing QR decomposition on $H$ takes $\mathcal{O}((i+1)ml)$. Note that we discard matrix $R$, only $Q$ is saved.
	\item Forming $T = A^T Q \in \mathbb{R}^{n\times (i+1)l}$ takes $\mathcal{O}((i+1)nl)$.
	\item Performing SVD on $T=\tilde{V}\tilde{\Sigma}W^T$ takes $\mathcal{O}((i+1)nl + 2[(i+1)l]^2)$, for matrix $\tilde{V}\in\mathbb{R}^{n\times (i+1)l}$ takes $\mathcal{O}((i+1)nl)$ and $\tilde{\Sigma}\in\mathbb{R}^{(i+1)l\times (i+1)l}$ takes $\mathcal{O}([(i+1)l]^2)$, and $W\in\mathbb{R}^{(i+1)l\times (i+1)l}$ takes $\mathcal{O}([(i+1)l]^2)$. In total, it takes $\mathcal{O}((i+1)nl + 2[(i+1)l]^2)$.
	\item Forming $\tilde{U}=QW$ takes $\mathcal{O}((i+1)ml)$.
\end{enumerate}
In total with $A$ taking $\mathcal{O}(mn)$, the space complexity of Krylov method is
\begin{align}
\mathcal{O}(mn+(3i+4)ml + (2i+3)nl + 2[(i+1)l]^2)
\end{align}

In practice, $l = 0.5n, i = 1$, then the space complexity of Krylov method will be
\begin{align}
\mathcal{O}(\dfrac{9}{2}mn + \dfrac{9}{2}n^2)
\end{align}

\subsubsection{Randomized PCA}
\begin{enumerate}
	\item Forming standard normal distribution matrix $G$ of size $n \times l$ takes $\mathcal{O}(nl)$, where in practice $l=0.5n$.
	\item Forming matrix $H^{(0)}=AG \in \mathbb{R}^{m\times l}$ takes $\mathcal{O}(ml)$.
	\item Forming matrix $H^{(i)}=A(A^T H^{(i-1)})\in \mathbb{R}^{m\times l}$ takes $\mathcal{O}(ml)$.
	\item Performing QR decomposition on $H$ takes $\mathcal{O}(ml)$. Note that we discard matrix $R$, only $Q$ is saved.
	\item Forming $T = A^T Q \in \mathbb{R}^{n\times l}$ takes $\mathcal{O}(nl)$.
	\item Performing SVD on $T=\tilde{V}\tilde{\Sigma}W^T$ takes $\mathcal{O}(nl + 2l^2)$, for matrix $\tilde{V}\in\mathbb{R}^{n\times l}$ takes $\mathcal{O}(nl)$ and $\tilde{\Sigma}\in\mathbb{R}^{l\times l}$ takes $\mathcal{O}(l^2)$, and $W\in\mathbb{R}^{l\times l}$ takes $\mathcal{O}(l^2)$. In total, it takes $\mathcal{O}(nl + 2l^2)$.
	\item Forming $\tilde{U}=QW$ takes $\mathcal{O}(ml)$.
\end{enumerate}
In total with $A$ taking $\mathcal{O}(mn)$, the space complexity of Randomized PCA is
\begin{align}
\mathcal{O}(mn+3ml+3nl+2l^2)
\end{align}

In practice, $l = 0.5n, i = 1$, then the space complexity of Randomized PCA will be
\begin{align}
\mathcal{O}(\dfrac{9}{2}mn + 3n^2)
\end{align}

\subsubsection{Truncated SVD}
We discuss the space complexity of Truncated SVD for each step.
\begin{enumerate}
	\item Forming matrix $A^T A\in\mathbb{R}^{n\times n}$ takes $\mathcal{O}(n^2)$.
	\item Performing eigen-decomposition on $A^TA \in\mathbb{R}^{n\times n}$ takes $\mathcal{O}(n^2+n)$.
	\item Taking the square root of each eigenvalue of $A^T A$ takes $\mathcal{O}(n)$.
	\item Forming $u_i = \dfrac{Av_i}{\sigma_i}$ takes $\mathcal{O}(mn)$, as each $u_i$ takes $\mathcal{O}(m)$ and we have $n$ equations like this, thus in total it takes $\mathcal{O}(mn)$.
	\item Forming $V$ takes $\mathcal{O}(n^2)$.
	\item Storing $n$ singular values takes $\mathcal{O}(n)$.
\end{enumerate}

In total with $A$ taking $\mathcal{O}(mn)$, the space complexity of Truncated SVD is
\begin{align}
\mathcal{O}(3n^2+3n+2mn)
\end{align}

\subsection{Summary of Complexity Analysis}

\begin{table}[h]
	\caption{Comparison of complexity}
	\label{table-complexity}
	\centering
	\begin{tabular}{lll}
		\toprule
		Method     & Time complexity   & Space complexity  \\
		\midrule
		Krylov method & 
		$\mathcal{O}(\dfrac{n^2}{2}+\dfrac{5}{2}mn^2+4 (m^2 n^2 +2mn^2+n^4 -\dfrac{2}{3}n^3))$  & $\mathcal{O}(\dfrac{9}{2}mn + \dfrac{9}{2}n^2)$ \\
		&  &  \\
		Randomized PCA & $\mathcal{O}(\dfrac{n^2}{2}+\dfrac{11}{4}mn^2+\dfrac{n^3}{6})$ & $\mathcal{O}(\dfrac{9}{2}mn + 3n^2)$       \\
		& & \\
		Truncated SVD & $\mathcal{O}(2mn^2+n^3+n+mn)$						& $\mathcal{O}(3n^2+3n+2mn)$ \\
		\bottomrule
	\end{tabular}
\end{table}

We summarized the time complexity and space complexity in Table \ref{table-complexity}.

Under the assumptions that $m \gg n$, for time complexity, by keeping the highest order term and its coefficient, we could see that for Krylov method, it takes $\mathcal{O}(4m^2n^2)$ FLOP while Randomized PCA takes $\mathcal{O}(\dfrac{11}{4}mn^2)$ FLOP and Truncated SVD only takes $\mathcal{O}(2mn^2)$ FLOP. Therefore, Truncated SVD is the fastest SVD algorithm among the aforementioned. Furthermore, Truncated SVD keeps all the eigenpairs rather than only first $k$ pairs as Krylov method and Randomized PCA do. 

For space complexity, we could see that Truncated SVD needs the least memory usage as $\mathcal{O}(3n^2+3n+2mn)$ and Krylov method needs the most memory space as $\mathcal{O}(\dfrac{9}{2}mn + \dfrac{9}{2}n^2)$. Randomized PCA holds the space complexity in between.

\section{Experiments}
We generate a matrix $A$ whose entry obeys standard normal distribution, i.e., $a_{ij} \sim \mathcal{N}(0,1)$, with 5 row sizes in list $[2000, 4000, 6000, 8000, 10000]$ and  12 column sizes in $[100, 200, 300, 400, 500, 600, 800, 900, 1000, 1200, 1500, 2000]$, 60 matrices of $A$ in total. The experiment is repeated 10 times to get an average runtime. On evaluating the residual error, the rate of Frobenius norm is used
\begin{align}
\delta = \dfrac{\|A-U\Sigma V^T \|_F}{\|A\|_F}
\end{align}

\begin{figure}
	\centering
	\includegraphics[width=0.65\linewidth]{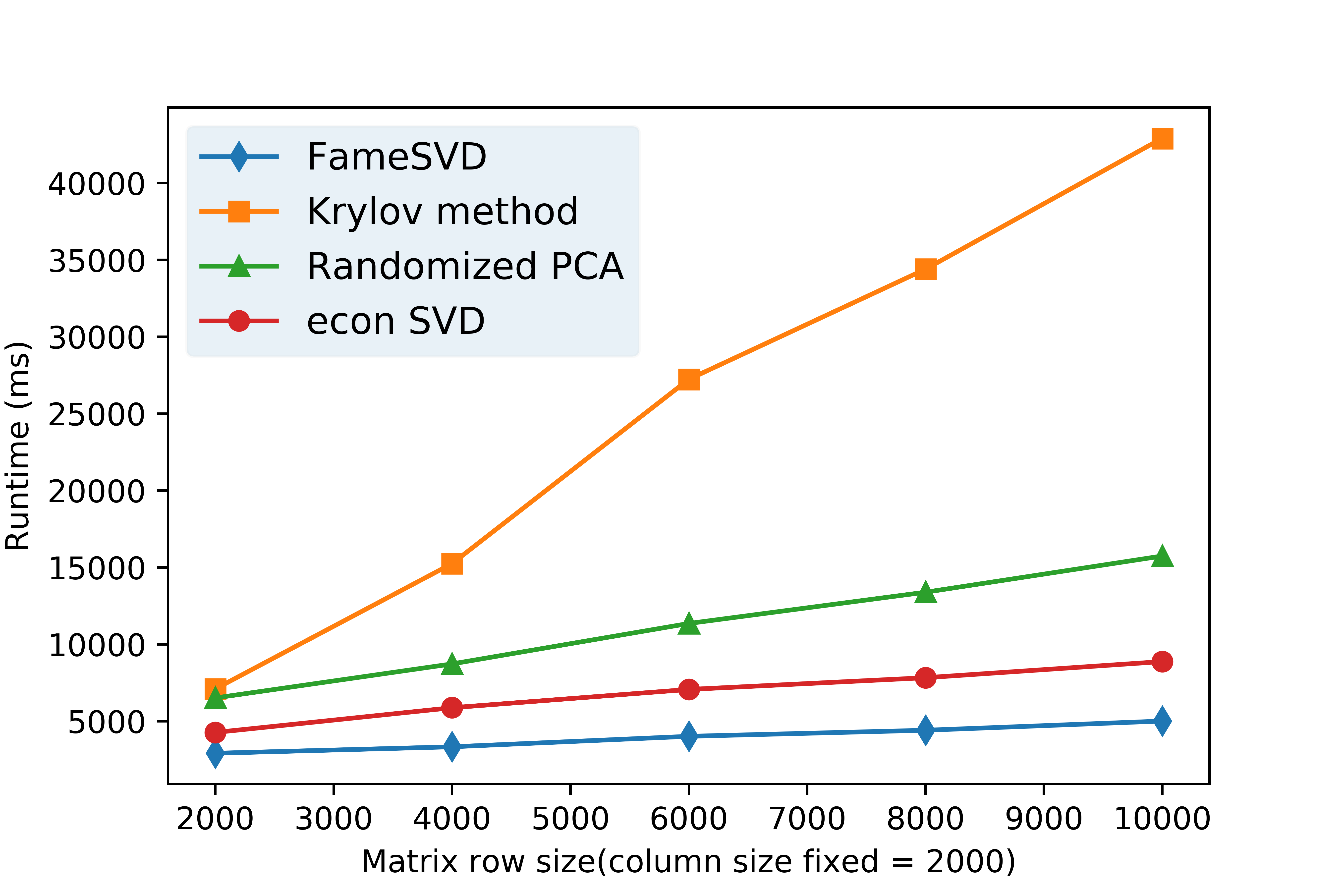}
	\caption{Runtime comparisons between 4 methods of SVD. The matrix column size is fixed at 2000, but row size varies. Truncated SVD method outperforms the rest.}
	\label{runtimefixcol}
\end{figure}

In Fig. \ref{runtimefixcol}, we compare the runtime of 4 SVD methods: Truncated SVD (FameSVD), Krylov method, Randomized PCA, econ SVD (MatLab built-in economic SVD). The matrix column size is fixed at 2000, and we increase the row size gradually. We could observe that all 4 methods follow a linear runtime pattern when row size increases. Of these 4 methods, Truncated SVD method outperforms the other 3 approaches.

\graphicspath{{E:/LXC/Graduate2Spring/FastSVD/experiment/MatLab/Fig/}}

\begin{figure}[hbt]
	\centering
	\includegraphics[width=0.65\linewidth]{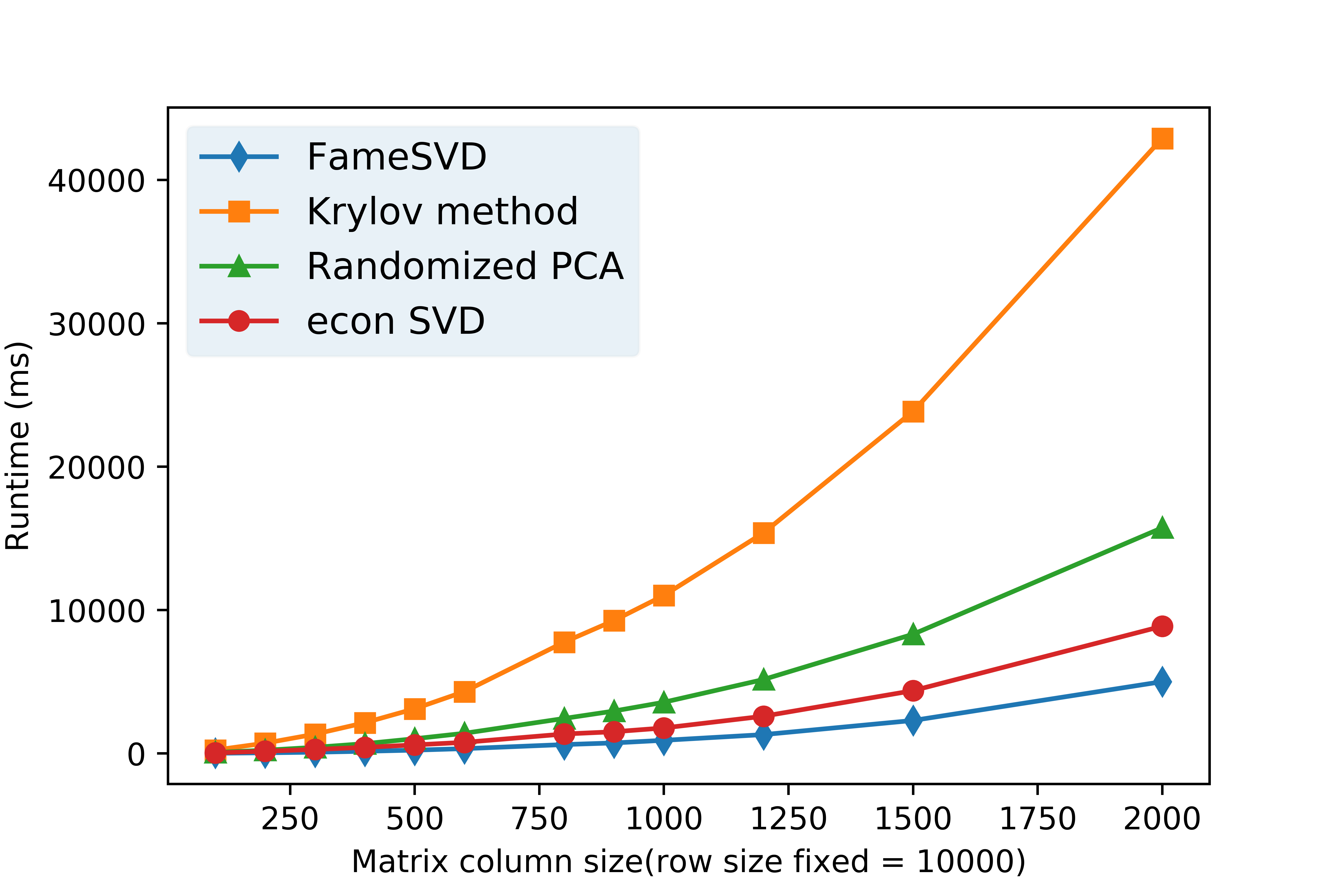}
	\caption{Runtime comparisons between 4 methods of SVD. The matrix row size is fixed at 10000, but column size varies. Truncated SVD method outperforms the rest.}
	\label{runtimefixrow}
\end{figure}
In Fig. \ref{runtimefixrow}, we fix the row size of matrix at 10000, and we increase the column size gradually. We could observe that all 4 methods behave as non-linear runtime pattern when row size increases. Out of all 4 methods, truncated SVD method takes the least runtime in every scenario.

\begin{figure}[!ht]
	\centering
	\includegraphics[width=0.65\linewidth]{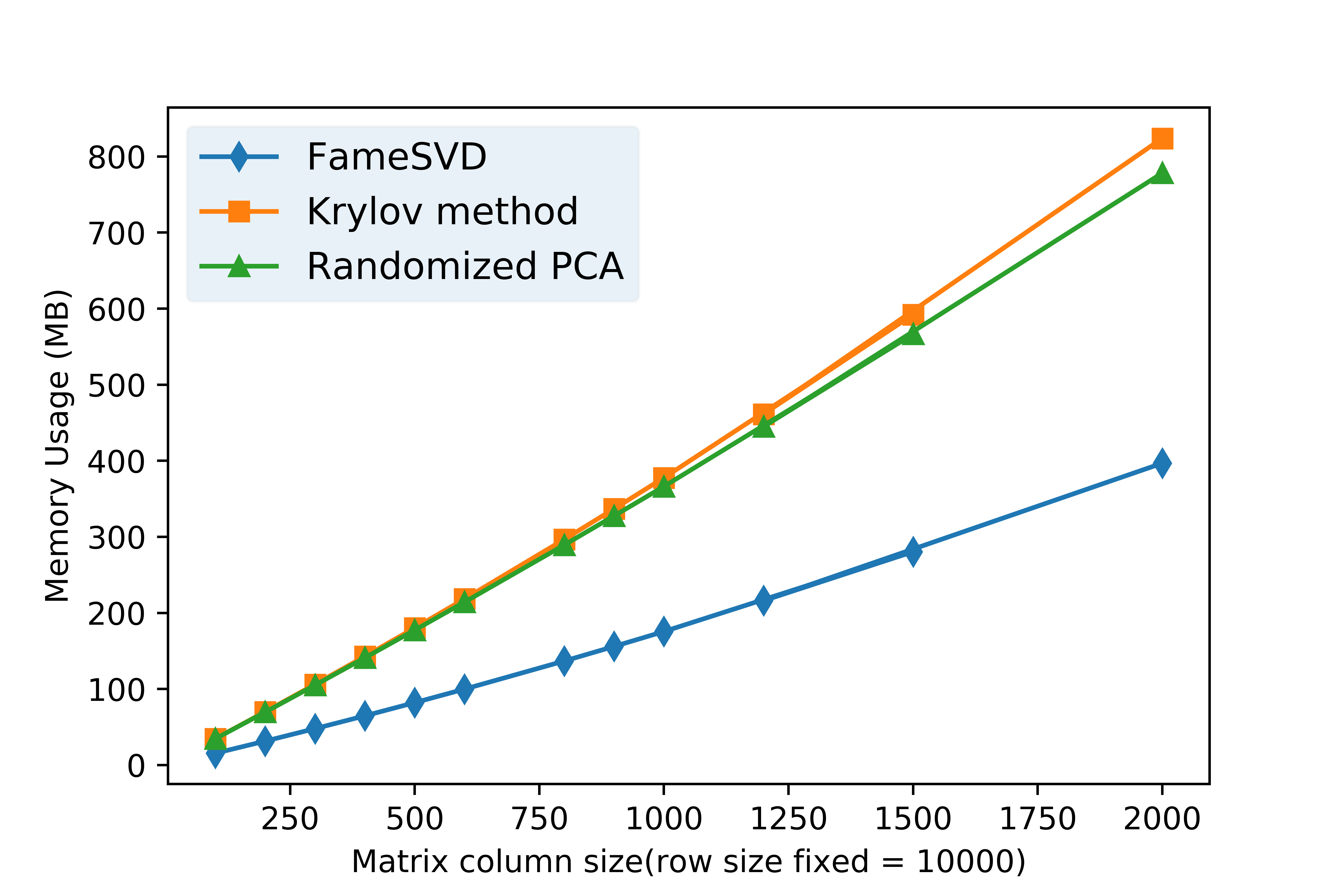}
	\caption{Memory usage between 3 methods of SVD. The matrix row size is fixed at 10000, but column size varies. Truncated SVD needs the least memory.}
	\label{memoryusagefig}
\end{figure}

In Fig. \ref{memoryusagefig}, the row size of matrix is fixed at 10000, but column size varies. We could see that truncated SVD uses the minimal amount of memory while Randomized PCA needs the most.

\begin{figure}[h]
	\centering
	\includegraphics[width=0.8\linewidth]{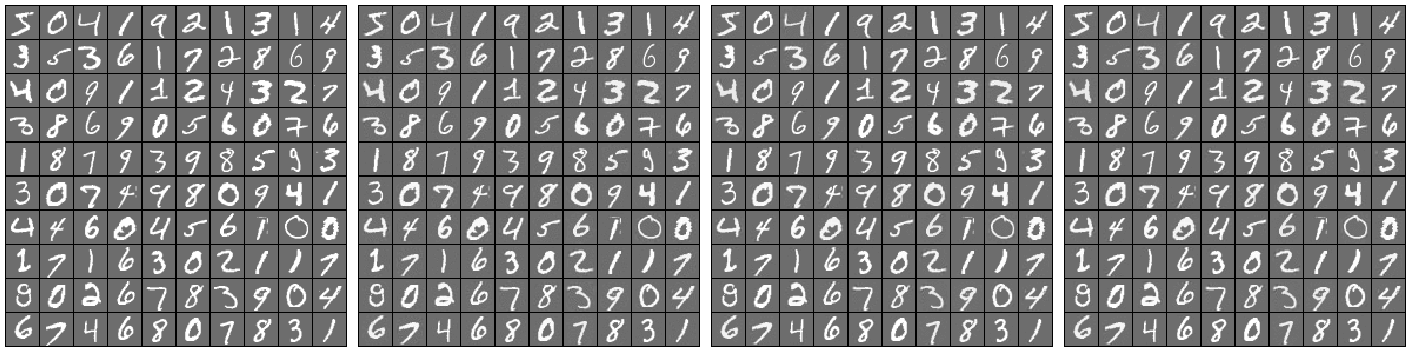}
	\caption{Left to right: Original image, Randomized PCA, Krylov method and truncated SVD. The first 392 (784/2) principal components are reserved.}
	\label{mnistfig}
\end{figure}

We also evaluate our algorithm on handwritten digit dataset MNIST\cite{Deng2012The}. We form our matrix $A$ with size $60000\times 784$ by concatenating $60000$ of vectorized $28\times 28$ intensity image. For runtime, it takes 4.54$s$ and 10.79$s$ for Randomized PCA and Krylov method respectively to obtain the first 392 (784/2) principal components. However, it only takes Truncated SVD \textbf{3.12s} to get all the 784 eigenvalues and eigenvectors; For memory usage, 1629.1MB for Randomized PCA and 1636.1MB for Krylov method. In the meanwhile, only \textbf{731.9MB} is used for truncated SVD. 

Our experiments are conducted on MatLab R2013a and Python 3.7 with NumPy 1.15.4, with Intel(R) Core(TM) i7-6700 CPU @ 3.40GHz 3.40GHz, 8.00GB RAM, and Windows 7. The truncated SVD is faster than the built-in economic SVD of both MatLab and NumPy.

\section{Conclusion}
the regular SVD performs the worst because it needs the most of memory usage and time. When all eigenpairs are needed ($k = n$), truncated SVD outperforms other SVD variants mentioned in this tutorial. The memory usage of these SVD methods grows linearly with matrix column size, but truncated SVD has the lowest growth rate. The runtime of truncated SVD grows sublinearly while Krylov method grows exponentially.

{\small
	\bibliographystyle{ieee}
	\bibliography{egbib}

\begin{thebibliography}{10}\itemsep=-1pt

\bibitem{Belabbas2009On}
M.~A. Belabbas and P.~J. Wolfe.
\newblock On sparse representations of linear operators and the approximation
  of matrix products.
\newblock 2009.

\bibitem{Belkin2014Laplacian}
M.~Belkin and P.~Niyogi.
\newblock Laplacian eigenmaps for dimensionality reduction and data
  representation.
\newblock {\em Neural Computation}, 15(6):1373--1396, 2014.

\bibitem{Cai2008A}
J.~F. Cai, E.~J. Candes, and Z.~Shen.
\newblock A singular value thresholding algorithm for matrix completion.
\newblock {\em Siam Journal on Optimization}, 20(4):1956--1982, 2008.

\bibitem{Cand2010The}
Cand, E.~J. S, and T.~Tao.
\newblock The power of convex relaxation: near-optimal matrix completion.
\newblock {\em IEEE Transactions on Information Theory}, 56(5):2053--2080,
  2010.

\bibitem{Candes2010Matrix}
E.~J. Candes and Y.~Plan.
\newblock Matrix completion with noise.
\newblock {\em Proceedings of the IEEE}, 98(6):925--936, 2010.

\bibitem{Chong2016Feature}
P.~Chong, K.~Zhao, Y.~Ming, and C.~Qiang.
\newblock Feature selection embedded subspace clustering.
\newblock {\em IEEE Signal Processing Letters}, 23(7):1018--1022, 2016.

\bibitem{deerwester1990indexing}
S.~Deerwester, S.~T. Dumais, G.~W. Furnas, T.~K. Landauer, and R.~A. Harshman.
\newblock Indexing by latent semantic analysis.
\newblock {\em Journal of the Association for Information Science and
  Technology}, 41(6):391--407, 1990.

\bibitem{Deng2012The}
L.~Deng.
\newblock The mnist database of handwritten digit images for machine learning
  research [best of the web].
\newblock {\em IEEE Signal Processing Magazine}, 29(6):141--142, 2012.

\bibitem{Frieze1998Fast}
A.~Frieze, R.~Kannan, and S.~Vempala.
\newblock Fast monte-carlo algorithms for finding low-rank approximations.
\newblock In {\em Symposium on Foundations of Computer Science}, 1998.

\bibitem{Golub1970Singular}
G.~H. Golub and C.~Reinsch.
\newblock Singular value decomposition and least squares solutions.
\newblock {\em Numerische Mathematik}, 14(5):403--420, 1970.

\bibitem{Halko2011An}
N.~Halko, P.~G. Martinsson, Y.~Shkolnisky, and M.~Tygert.
\newblock An algorithm for the principal component analysis of large data sets.
\newblock {\em Siam Journal on Scientific Computing}, 33(5):2580--2594, 2011.

\bibitem{Halko2010Finding}
N.~Halko, P.~G. Martinsson, and J.~A. Tropp.
\newblock Finding structure with randomness: Probabilistic algorithms for
  constructing approximate matrix decompositions.
\newblock {\em Siam Review}, 53(2):217--288, 2010.

\bibitem{hofmann2001unsupervised}
T.~Hofmann.
\newblock Unsupervised learning by probabilistic latent semantic analysis.
\newblock {\em Machine Learning}, 42(1):177--196, 2001.

\bibitem{J1989Estimating}
H.~W. J.~Kuczynski.
\newblock Estimating the largest eigenvalue by the power and lanczos algorithms
  with a random start.
\newblock {\em Siam Journal on Matrix Analysis \& Applications},
  13(4):1094--1122, 1989.

\bibitem{Kai2010Clustered}
Z.~Kai and J.~T. Kwok.
\newblock Clustered nyström method for large scale manifold learning and
  dimension reduction.
\newblock {\em IEEE Transactions on Neural Networks}, 21(10):1576, 2010.

\bibitem{Partridge2002Robust}
M.~Partridge and M.~Jabri.
\newblock Robust principal component analysis.
\newblock In {\em Neural Networks for Signal Processing X, IEEE Signal
  Processing Society Workshop}, 2002.

\bibitem{Recht2012Exact}
B.~Recht.
\newblock {\em Exact matrix completion via convex optimization}.
\newblock 2012.

\bibitem{Tipping2010Probabilistic}
M.~E. Tipping and C.~M. Bishop.
\newblock Probabilistic principal component analysis.
\newblock {\em Journal of the Royal Statistical Society}, 61(3):611--622, 2010.

\bibitem{Woodruff2014Low}
D.~P. Woodruff.
\newblock Low rank approximation lower bounds in row-update streams.
\newblock In {\em International Conference on Neural Information Processing
  Systems}, 2014.

\bibitem{Woolfe2008A}
F.~Woolfe, E.~Liberty, V.~Rokhlin, and M.~Tygert.
\newblock A fast randomized algorithm for the approximation of matrices.
\newblock {\em Applied \& Computational Harmonic Analysis}, 25(3):335--366,
  2008.

\bibitem{Xu1995Robust}
.~Xu, L. and A.~L. Yuille.
\newblock Robust principal component analysis by self-organizing rules based on
  statistical physics approach.
\newblock {\em IEEE Transactions on Neural Networks}, 6(1):131--43, 1995.

\bibitem{Zhang2010Regularized}
Z.~Zhang, G.~Dai, C.~Xu, and M.~I. Jordan.
\newblock Regularized discriminant analysis, ridge regression and beyond.
\newblock {\em Journal of Machine Learning Research}, 11(12):2199--2228, 2010.

\bibitem{Zou2006Sparse}
H.~Zou, T.~Hastie, and R.~Tibshirani.
\newblock Sparse principal component analysis.
\newblock {\em Journal of Computational and Graphical Statistics},
  15(2):265--286, 2006.

\end{thebibliography}
}
\end{document}